\documentclass[12pt]{amsart}         
\usepackage{amscd}                   
\usepackage{amssymb}                 

\newtheorem{thm}{Theorem}[section]
\newtheorem{lem}[thm]{Lemma}

\renewcommand{\thestep}{}
\newtheorem{prop}[thm]{Proposition}

\theoremstyle{definition}

\renewcommand{\thecase}{}

\newtheorem{rmk}[thm]{Remark}

\theoremstyle{remark}

\makeatletter
\def\alphenumi{
  \def\theenumi{\alph{enumi}}
  \def\p@enumi{\theenumi}
  \def\labelenumi{(\@alph\c@enumi)}}
\makeatother

\makeatletter
\def\thecase{\@arabic\c@case}
\makeatother

\numberwithin{equation}{section}

\makeatletter
\def\thestep{\@arabic\c@step}
\makeatother

\newenvironment{pf*}[1]{\begin{proof}[#1]}{\end{proof}}

\newcommand\barB{{\bar B}}

\newcommand\barM{{\bar{M}}}

\newcommand\fs{{\mathfrak{S}}}
\newcommand\fS{{\underline{\mathfrak{S}}}}

\newcommand\ubarfV{{\underline{\mathfrak{V}}}}

\newcommand\AAA{\mathbb{A}}

\newcommand\CC{\mathbb{C}}

\newcommand\RR{\mathbb{R}}

\newcommand\ZZ{\mathbb{Z}}

\newcommand\bM{{\mathbf{M}}}

\newcommand\bx{{\mathbf{x}}}

\newcommand\ssG{{{}^\circ\mathcal{G}}}

\newcommand{\cov}{\nabla}

\newcommand\half{{\textstyle{\frac{1}{2}}}}

\newcommand\quarter{{\textstyle{\frac{1}{4}}}}
\newcommand\threehalf{{\textstyle{\frac{3}{2}}}}

\newcommand\sixtyfourth{{\textstyle{\frac{1}{{64}}}}}

\newcommand\ff{{\mathfrak{f}}}

\newcommand\fh{{\mathfrak{h}}}

\newcommand\fm{{\mathfrak{m}}}
\newcommand\fM{{\mathfrak{M}}}

\newcommand\fV{{\mathfrak{V}}}

\newcommand\eps{\varepsilon}
\newcommand\ga{\gamma}

\newcommand\la{\lambda}
\newcommand\La{\Lambda}

\newcommand\Om{\Omega}

\newcommand\gl{{\mathfrak{g}\mathfrak{l}}}
\newcommand\fsl{{\mathfrak{s}\mathfrak{l}}}
\newcommand\hol{{\mathfrak{h}\mathfrak{o}\mathfrak{l}}}
\newcommand\so{{\mathfrak{s}\mathfrak{o}}}

\newcommand\su{{\mathfrak{s}\mathfrak{u}}}
\newcommand\fu{{\mathfrak{u}}}

\newcommand\PU{\operatorname{PU}}

\newcommand\SO{\operatorname{SO}}

\newcommand\SU{\operatorname{SU}}
\newcommand\U{\operatorname{U}}

\newcommand\less{\setminus}

\newcommand{\8}{\infty}

\newcommand\ad{{\operatorname{ad}}}

\newcommand\asd{{\operatorname{asd}}}

\newcommand\dist{\operatorname{dist}}

\newcommand\End{\operatorname{End}}

\newcommand\Hol{\operatorname{Hol}}
\newcommand\Hom{\operatorname{Hom}}

\newcommand\loc{\operatorname{loc}}

\newcommand\red{\operatorname{red}}

\newcommand\Sym{\operatorname{Sym}}

\newcommand\id{{\mathrm{id}}}

\newcommand\spinc{\text{$\text{spin}^c$ }}

\newcommand\sA{{\mathcal{A}}}
\newcommand\sB{{\mathcal{B}}}
\newcommand\sC{{\mathcal{C}}}

\newcommand\sG{{\mathcal{G}}}

\newcommand\sP{{\mathcal{P}}}

\newcommand\tsC{{\tilde\sC}}

\newcommand\vecfm{{\vec{\fm}}}
\newcommand\vectau{{\vec{\tau}}}
\newcommand\vecvartheta{{\vec{\vartheta}}}

\begin{document}
\title[Uhlenbeck Compactness and Transversality for $PU(2)$ Monopoles]
{Uhlenbeck Compactness and Transversality for the Moduli Space of
PU(2) Monopoles$\text{}^1$}  
\author[Paul M. N. Feehan]{Paul M. N. Feehan}
\address{Department of Mathematics\\
Harvard University\\
One Oxford Street\\
Cambridge, MA 02138}
\email{feehan@math.harvard.edu}
\author[Thomas G. Leness]{Thomas G. Leness}
\address{Department of Mathematics\\
Michigan State University\\
East Lansing, MI 48824--1027}
\email{leness@math.msu.edu}
\thanks{The first author was supported in part by an NSF Mathematical 
Sciences Postdoctoral Fellowship under grant DMS 9306061.}
\footnotetext[1]{Submitted to a print journal,
June 2, 1997. This version: October 12, 1997. \texttt{dg-ga/9710012}}
\maketitle


\section{Introduction}
This research announcement reports on some results in \cite{FL1} concerning
the existence of perturbations for the $\PU(2)$ monopole equations,
yielding both useful transversality properties and an Uhlenbeck
compactification for this perturbed moduli space. 

A method has been proposed in \cite{OTQuaternion,PTCambridge,PTLocal} to
prove Witten's conjecture concerning the relation between the Donaldson and
Seiberg-Witten invariants of smooth four-manifolds
\cite{DonSW,Witten,MooreWitten}.  The idea is to use a moduli space of 
solutions to the $\PU(2)$ monopole equations, which are a natural
generalization of the $\U(1)$ monopole equations of Seiberg and Witten and the
anti-self-dual equation for $\SO(3)$ connections, to construct a cobordism
between links of compact moduli spaces of $\U(1)$ monopoles of
Seiberg-Witten type and the moduli space of anti-self-dual
connections, which appear as singularities in this larger moduli space.
The conjecture holds for all four-manifolds whose Donaldson and
Seiberg-Witten invariants have been independently computed by direct
calculation. A basic requirement of this cobordism technique is the
existence of an Uhlenbeck compactification for the moduli space of $\PU(2)$
monopoles and of generic-parameter transversality results for all the
moduli spaces of $\PU(2)$ monopoles which appear in this compactification,
at least away from the anti-self-dual and $\U(1)$ solutions.

In \S \ref{sec:HolonomyPerturbations} we describe the holonomy
perturbations we use in order to achieve transversality for the moduli
space of solutions to the perturbed $\PU(2)$ monopole equations. In \S
\ref{sec:UhlenbeckCompactness} and \S \ref{sec:Transversality} we
state our Uhlenbeck compactness and transversality theorems for this moduli
space and outline the proofs from \cite{FL1}, to which we refer for
detailed arguments.

\section{Holonomy perturbations}
\label{sec:HolonomyPerturbations}
We consider Hermitian two-plane bundles $E$ over $X$ whose determinant line
bundles $\det E$ are isomorphic to a fixed Hermitian line bundle over $X$
endowed with a fixed $C^\8$, unitary connection $A_e$.  Let
$(\rho,W)$ be a \spinc structure on $X$, where
$\rho:T^*X\to\End W$ is the Clifford map, and the Hermitian
four-plane bundle $W=W^+\oplus W^-$ is endowed with a $C^\8$ \spinc
connection. {\em The \spinc structure $(\rho,W)$, the \spinc
connection on $W$, and the Hermitian line bundle together with
its connection  $A_e$ are fixed once and for all.\/}

Let $k\ge 2$ be an integer and let $\sA_E$ be the space of $L^2_k$
connections $A$ on the $\U(2)$ bundle $E$ all inducing the fixed
determinant connection on $\det E$.  Equivalently, following \cite[\S
2(i)]{KMStructure}, we may view $\sA_E$ as the space of $L^2_k$ connections
$A$ on the $\PU(2)=\SO(3)$ bundle $\su(E)$.  We shall often pass back and
forth between these viewpoints, via the fixed connection on $\det E$,
relying on the context to make the distinction clear.  Let
$D_A:L^2_k(W^+\otimes E)\to L^2_{k-1}(W^-\otimes E)$ be the corresponding
Dirac operators.  Given a connection $A$ on $E$ with curvature
$F_A\in L^2_{k-1}(\La^2\otimes\fu(E))$, then $(F_A^+)_0 \in
L^2_{k-1}(\La^+\otimes\su(E))$ denotes the traceless part of its self-dual
component. Equivalently, if $A$ is a connection on $\su(E)$ with
curvature $F_A\in L^2_{k-1}(\La^2\otimes\so(\su(E)))$, then
$\ad^{-1}(F_A^+) \in L^2_{k-1}(\La^+\otimes\su(E))$ is its self-dual
component, viewed as a section of $\La^+\otimes\su(E)$ via the isomorphism
$\ad:\su(E)\to\so(\su(E))$. When no confusion can arise, the
isomorphism $\ad:\su(E)\to\so(\su(E))$ will be implicit and so we regard
$F_A$ as a section of $\La^+\otimes\su(E)$ when $A$ is a connection on
$\su(E)$.

For an $L^2_k$ section $\Phi$ of $W^+\otimes E$, let $\Phi^*$ be its
pointwise Hermitian dual and let $(\Phi\otimes\Phi^*)_{00}$ be the
component of the Hermitian endomorphism $\Phi\otimes\Phi^*$ of $W^+\otimes
E$ which lies in $\su(W^+)\otimes\su(E)$. The Clifford multiplication
$\rho$ defines an isomorphism $\rho^+:\La^+\to\su(W^+)$ and thus an
isomorphism $\rho^+=\rho^+\otimes\id_{\su(E)}$ of $\La^+\otimes\su(E)$ with
$\su(W^+)\otimes\su(E)$. Then
\begin{align}
F_A^+ - (\rho^+)^{-1}(\Phi\otimes\Phi^*)_{00} &= 0, 
\label{eq:IntroUnpertPT}\\
D_A\Phi &= 0, \notag
\end{align}
are the unperturbed equations considered in
\cite{OTVortex,OTQuaternion,PTCambridge,PTLocal} (the trace conditions
vary slightly --- see \cite{FL1}), for a pair $(A,\Phi)$ consisting of a
connection on $\su(E)$ and a section $\Phi$ of $W^+\otimes
E$. Equivalently, given a pair $(A,\Phi)$ with fixed-determinant connection
$A$ on $E$, the equations
\eqref{eq:IntroUnpertPT} take the same form except that $F_A^+$ is replaced
by $(F_A^+)_0$.

In this section we briefly describe our
holonomy perturbations of these equations 
\cite{FL1}. These perturbations allow us to prove transversality for the
moduli space of solutions, away from points where the connection is
reducible or the spinor vanishes identically,
and to prove the existence of an Uhlenbeck
compactification for this perturbed moduli space.

Donaldson's proof of the connected-sum theorem for his polynomial
invariants \cite[Theorem B]{DonPoly} makes use of certain `extended
anti-self-dual equations' \cite[Equation (4.24)]{DonPoly}
to which the Freed-Uhlenbeck generic metrics
theorem does not apply \cite[\S 4(v)]{DonPoly}. These extended equations
model a neighborhood of the product connection appearing in the Uhlenbeck
compactification of the moduli space of anti-self-dual $\SU(2)$ connections.
To obtain transversality for the zero locus of
the extended equations, he employs holonomy perturbations which give
gauge equivariant $C^\8$ maps $\sA_E^*\to
\Om^+(\su(E))$ \cite[\S 2]{DonOrient}, \cite[pp. 282--287]{DonPoly}. 
These perturbations are continuous with respect to
Uhlenbeck limits and yield transversality not only for the top stratum,
but also for all lower strata and for all intersections of the geometric
representatives defining the Donaldson invariants. 

In \cite[\S 2.4 \& Appendix]{FL1} we describe a generalization of
Donaldson's idea which we use to prove transversality for the moduli space
of solutions to a perturbed version of the $\PU(2)$ monopole equations
\eqref{eq:IntroUnpertPT}.
Unfortunately, in the case of the moduli space of $\PU(2)$ monopoles, the
analysis is considerably more intricate. In Donaldson's application, some
important features ensure that the requisite analysis is relatively
tractable: (i) reducible connections can be excluded from the
compactification of the extended moduli spaces \cite[p. 283]{DonPoly}, (ii)
the cohomology groups for the elliptic complex of his extended equations
have simple weak semi-continuity properties with respect to Uhlenbeck
limits
\cite[Proposition 4.33]{DonPoly}, and (iii) the zero locus being perturbed
is cut out of a finite-dimensional manifold \cite[p. 281, Lemma 4.35, \&
Corollary 4.38]{DonPoly}. For the development of Donaldson's method for
$\PU(2)$ monopoles described here and in detail in \cite{FL1}, none of
these simplifying features hold and so the corresponding transversality
argument is rather complicated. Indeed, one can see from
Proposition 7.1.32 in
\cite{DK} that because of the Dirac operator, the behavior of the cokernels
of the linearization of the $\PU(2)$ monopole equations can be quite
involved under Uhlenbeck limits. The method we describe below uses an
infinite sequence of perturbing sections defined on the
infinite-dimensional configuration space of pairs; when restricted to small
enough open balls in the configuration space, away from reducibles, only
finitely many of these perturbing sections are non-zero and they vanish
along the reducibles.

Let $\sG_E$ be the Hilbert Lie group of $L^2_{k+1}$ unitary gauge
transformations of $E$ with {\em determinant one\/}. It is often
convenient to take quotients by a slightly larger symmetry group than
$\sG_E$ when discussing pairs, so let $S_Z^1$ denote the center of $\U(2)$
and set
$\ssG_E := S_Z^1\times_{\{\pm\id_E\}}\sG_E$, 
which we may view as the group of $L^2_{k+1}$ unitary gauge transformations
of $E$ with constant determinant.  The stabilizer of a unitary
connection on $E$ in $\ssG_E$ always contains the center
$S^1_Z\subset\U(2)$. We call $A$ {\em irreducible\/} if its stabilizer is
exactly $S_Z^1$ and {\em reducible\/} otherwise. Let $\sB_E(X) :=
\sA_E(X)/\sG_E$ be the quotient space of $L^2_k$
connections on $E$ with fixed determinant connection and let $\sA^*_E(X)$
and $\sB_E^*(X)$ be the subspace of irreducible $L^2_k$ connections
and its quotient. As before, we may equivalently view $\sB_E(X)$ and
$\sB_E^*(X)$ as quotients of the spaces of $L^2_k$ connections on $\su(E)$
by the induced action of $\sG_E$ on $\su(E)$.

We fix $r\ge k+1$ and define gauge-equivariant $C^\8$ maps, whose
construction we outline below,
\begin{align}
&\sA_E(X) \ni A\mapsto \vectau\cdot\vecfm(A)
\in L^2_{k+1}(X,\gl(\La^+)\otimes_\RR\so(\su(E))), 
\label{eq:GaugeEquivariantMap}\\
&\sA_E(X) \ni A\mapsto \vecvartheta\cdot\vecfm(A)
\in L^2_{k+1}(X,\Hom(W^+,W^-)\otimes_\CC\fsl(E)), \notag
\end{align}
where $\vectau := (\tau_{j,l,\alpha})$ is a sequence in $C^r(X,\gl(\La^+))$
and $\vecvartheta := (\vartheta_{j,l,\alpha})$ is a sequence in
$C^r(X,\La^1\otimes\CC)$, while $\vecfm(A) := (\fm_{j,l,\alpha}(A))$ is a
sequence in $L^2_{k+1}(X,\su(E))$ of holonomy sections, and 
\begin{align*}
\vectau\cdot\vecfm(A)
&:= 
\sum_{j,l,\alpha}\tau_{j,l,\alpha}\otimes_\RR
\ad(\fm_{j,l,\alpha}(A)), \\
\vecvartheta\cdot\vecfm(A)
&:=
\sum_{j,l,\alpha}\rho(\vartheta_{j,l,\alpha})
\otimes_\CC\fm_{j,l,\alpha}(A).
\end{align*}
We digress briefly to outline the construction of the gauge equivariant
maps $\fm_{j,l,\alpha}$.

Let $\gamma\subset X$ be a $C^\8$ loop based
at a point $x_0\in X$ and let $\Hol_{\ga,x_0}(A)\in\SO(\su(E))|_{x_0}$ be the
holonomy of a smooth $\SO(3)$ connection $A$ around the loop $\gamma$. The
exponential map $\exp:\so(3)\to\SO(3)$ gives a diffeomorphism from a ball $2D$
around the origin in $\so(3)$ to a ball around the identity in $\SO(3)$. Let
$\psi:\RR\to[0,1]$ be a $C^\8$ cutoff function such that $\psi(|\zeta|)=1$
for $\zeta\in \half D$, $\psi(|\zeta|)>0$ for $\zeta\in D$, and
$\psi(|\zeta|)=0$ for $\zeta\in
\so(3)-D$ and define 
$$
\hol_{\gamma,x_0}(A) 
:= 
\psi\left(|\exp^{-1}(\Hol_{\ga,x_0}(A))|\right)
\cdot\ad^{-1}\exp^{-1}(\Hol_{\ga,x_0}(A)).
$$
Suppose $Y\subset X$ is a simply-connected open subset containing $x_0$.
If $A|_Y$ is irreducible then there are three loops $\{\gamma_l\}_{l=1}^3
\subset Y$, depending on $A$, such that the set
$\{\hol_{\gamma_l,x_0}(A)\}_{l=1}^3$ is a basis for 
$\su(E)|_{x_0} = \ad^{-1}(\so(\su(E)|_{x_0}))$.
We extend $\hol_{\gamma,x_0}(A)$ to a $C^\8$
section $\hat\fh_{\gamma}(A)$ of $\su(E)$ by radial parallel translation,
with respect to $A$ over a small ball $2B$ and then multiplying by
a $C^\8$ cutoff function $\varphi$ on $X$ which is positive on $B$
and identically zero on $X-B$. The 
set $\{\varphi\hat\fh_{\gamma_l}(A)|_y\}_{l=1}^3$ then spans
$\su(E)|_y$ for $y\in B$. 

For an $L^2_k$ unitary connection $A$ on $E$ with $k\ge 2$, the section
$\hat\fh_{\gamma}(A)$ need not be in $L^2_{k+1}$. So we use the Neumann
heat operator , for fixed small $t>0$,
to construct an $L^2_{k+1}$ section
$$
\fh_{\gamma}(A) := \exp(-td_{A|2B}^*d_{A|2B})\hat\fh_{\gamma}(A)
$$
of $\su(E)$ over $2B$ which converges to $\hat\fh_{\gamma}(A)$ in
$C^0(B)$ as $t\to 0$. Therefore, for a small enough time
$t(A)$, the set $\{\varphi\fh_{\gamma_l}(A)|_y\}_{l=1}^3$,
spans $\su(E)|_y$ for all $y\in B$ just as before.
If $A_0$ is an $L^2_k$ unitary connection on
$E|_{2B}$ and $\ga_1,\ga_2,\ga_3$ are loops in $2B$ based
at $x_0$ such that $\{\fh_{\gamma_l,x_0}(A_0)\}_{l=1}^3$ spans
$\su(E)|_{x_0}$, then there is a positive constant
$\eps(A_0,\{\ga_l\})$ such that if $A$
is an $L^2_k$ unitary connection on $E|_{2B}$ with
$\|A-A_0\|_{L^2_{k,A_0}(2B)} < \eps$,
then the set $\{\fh_{\gamma_l,x_0}(A)\}_{l=1}^3$ spans $\su(E)|_B$.

Given this digression, we can now construct the maps $\fm_{j,l,\alpha}$.
Let $\{4B_j\}_{j=1}^{N_b}$, be a disjoint collection of open balls with
centers $x_j$ and radius $R_0$, where $N_b$ is a fixed integer to be
determined later.  The quotient space $\sB^*_E(2B_j)$ of irreducible
$L^2_k$ unitary connections on $E|_{2B_j}$ is a paracompact $C^\8$ manifold
modelled on a separable Hilbert space. For each $j=1,\dots,N_b$ and each
point $[A_0]$ in $\sB_E^*(2B_j)$, we can find loops
$\{\ga_{j,l,A_0}\}_{l=1}^3$, contained in $2B_j$ and based at $x_j$ such
that $\{\fh_{\gamma_{j,l,A_0}}(A_0)\}_{l=1}^3$, spans $\fu(E)|_{B_j}$. For
each such point $[A_0]$, there is an $L^2_{k,A_0}$ ball $B_{[A_0]}(\eps_{A_0})
\subset \sB_E^*(2B_j)$ such that for all $[A]\in B_{[A_0]}(\eps_{A_0})$,
the set $\{\fh_{\gamma_{j,l,A_0}}(A)\}_{l=1}^3$ spans $\su(E)|_{B_j}$.
These balls give an open cover of $\sB^*_E(2B_j)$ and hence there is a
locally finite refinement of this open cover,
$\{U_{j,\alpha}\}_{\alpha=1}^\infty$, and a collection of suitably chosen
$C^\8$ cutoff functions $\chi_{j,\alpha}$ with supports in the $L^2_k$ balls
containing the $U_{j,\alpha}$ such
that
$$
\sum_{\alpha}\chi_{j,\alpha}[A]>0, \qquad
[A]\in \sB^*_E(2B_j), \quad j=1,\dots,N_b.
$$
Hence, for each $U_{j,\alpha}$, we obtain loops
$\{\gamma_{j,l,\alpha}\}_{l=1}^3\subset 2B_j$ such that for all $[A]\in
U_{j,\alpha}$, the sections $\fh_{\gamma_{j,l,\alpha}}(A)$ span
$\su(E)|_{B_j}$. The cutoff functions $\chi_{j,\alpha}$ are chosen so that
all their derivatives are bounded (in the obvious sense analogous to
\eqref{eq:NormDefn} below). 

Let $\beta$ be a $C^\8$ cutoff function on $\RR$ such that $\beta(t) = 1$
for $t\le \half$ and $\beta(t)=0$ for $t\ge 1$, with $\beta(t)>0$ for
$t<1$. Then the $C^\8$ gauge-invariant maps
$\sA_E(X)\to \RR$, $A\mapsto \beta_j[A]$, given by
\begin{equation}
\beta_j[A]
:=
\beta\left(\frac{1}{\eps_0^2}\int_{B(x_j,4R_0)}
\beta\left(\frac{\dist(\cdot,x_j)}{4R_0}\right)|F_A|^2\,dV\right)
\label{eq:ConnEnergyCutoff}
\end{equation}
are zero when the energy of the connection $A$ is greater than or equal to
$\half\eps_0^2$ over a ball $2B_j$. Here, $\eps_0(g,A_{\det W^+},A_{\det
E})$ is a universal constant \cite{FL1}.  Finally, we define $C^\8$ cutoff
functions on $X$ by setting $\varphi_j := \beta(\dist_g(\cdot,x_j)/R_0)$,
so that $\varphi_j$ is positive on the ball $B_j$ and zero on its
complement in $X$.  We define gauge-equivariant $C^\8$ maps $\sA_E(X)\to
L^2_{k+1}(X,\su(E))$, $A\mapsto \fm_{j,l,\alpha}(A)$, by setting
\begin{equation}
\fm_{j,l,\alpha}(A)
:= \beta_j[A]\chi_{j,\alpha}\left([A|_{2B_j}]\right)
\varphi_j\fh_{\gamma_{j,l,\alpha}}(A). 
\label{eq:HolonomySection} 
\end{equation}
Thus, at each point $A\in \sA_E(X)$ only a finite number of the
$\fm_{j,l,\alpha}(A)$ are non-zero and we show that each map
$\fm_{j,l,\alpha}$ is $C^\8$ with uniformly bounded derivatives of all
orders on $\sA_E(X)$ in the sense of \eqref{eq:NormDefn}. The energy cutoff
functions $\beta_j$ ensure continuity across the Uhlenbeck boundary.  The
universal bound on $\|F_A\|_{L^2(X)}^2$ in \eqref{eq:EnergyBound} below for
solutions $(A,\Phi)$ to the perturbed $\PU(2)$ monopole equations
\eqref{eq:PT} ensures that the number $N_b$ of
balls $B_j$ may be chosen sufficiently large that for every $\PU(2)$
monopole $(A,\Phi)$, there is at least one ball $B_{j'}$
whose associated holonomy sections $\{\fm_{j',l,\alpha}(A)\}_{l=1}^3$
span $\su(E)|_{B_{j'}}$. 

The parameters $\vectau$ and $\vecvartheta$ vary in the Banach spaces of
$\ell^1_\delta(\AAA)$ sequences in $C^r(X,\gl(\La^+))$ and
$C^r(X,\La^1\otimes\CC)$, 
respectively, where $\AAA = \{(j,l,\alpha)\}$, with norm
$$
\|\vecvartheta\|_{\ell^1_\delta(C^r(X))}
:= 
\sum_{j,l,\alpha}\delta_\alpha^{-1}
\|\vartheta_{j,l,\alpha}\|_{C^r(X)},
$$
and similarly for $\|\vectau\|_{\ell^1_\delta(C^r(X))}$.
For any open subset $U\subset \sA_E(X)$ and $C^s$ map $\ff:U\to
L^2_{k+1}(X,\su(E))$, we define
\begin{equation}
\|\ff\|_{C^s(U)}
:=
\sup_{A\in U}
\sup_{\substack{1\le i\le s \\ \|a_i\|_{L^2_{k,A}(X)}\le 1}}
\|(D^s\ff)_A(a_1,\dots,a_s)\|_{L^2_{k+1,A}(X)},
\label{eq:NormDefn}
\end{equation}
and similarly for maps to
$L^2_{k+1}(X,\gl(\La^+)\otimes\so(\su(E)))$ and
$L^2_{k+1}(X,\Hom(W^+,W^-)\otimes\fsl(E))$. The $C^s$ norm of maps on
$\sA_E(X)$ defined by \eqref{eq:NormDefn} is gauge invariant and so
descends to a $C^s$ norm on sections over the quotient $\sB_E^*(X)$.

\begin{prop}
\label{prop:SmoothAcrossReducibles}
\cite{FL1}
There exists a sequence $\delta:=(\delta_\alpha)_{\alpha=1}^\8 \in
\ell^\8((0,1])$ of positive weights
such that the gauge-equivariant maps in \eqref{eq:GaugeEquivariantMap}
are $C^\8$, with uniformly bounded differentials of all orders,
$$
\|\vecvartheta\cdot\vecfm\|_{C^s(\sA_E)}
\le C\|\vecvartheta\|_{\ell^1_\delta(C^r(X))}
\quad\text{and}\quad
\|\vectau\cdot\vecfm\|_{C^s(\sA_E)}
\le C\|\vectau\|_{\ell^1_\delta(C^r(X))},
$$
for all $s\ge 0$ and positive constants $C=C(g,k,s)$.
\end{prop}

The parameters $\tau_0,\vartheta_0,\vectau,\vecvartheta$ can be chosen so that
\begin{align}
\|\tau_0\|_{L^\8(X)}+
\sup_{A\in\sA_E}\|\vectau\cdot\vecfm(A)\|_{L^\8(X)} &< \sixtyfourth, 
\label{eq:CompactEstPertC0}
\\
\|\vartheta_0\|_{L^\8_1(X)}
\sup_{A\in\sA_E}\|\vecvartheta\cdot\vecfm(A)\|_{L^\8_{1,A}(X)} &\le 1,
\notag
\end{align}
{\em provided $k\ge 3$} and we shall therefore assume this constraint is in
effect for the remainder of the article. These bounds are used in
\S \ref{sec:UhlenbeckCompactness} and they
in turn follow from the $C^0$ estimates of
Proposition \ref{prop:SmoothAcrossReducibles} for
$$
\|\vectau\|_{\ell^1_\delta(C^r(X))} \le \eps_\tau
\quad\text{and}\quad
\|\vecvartheta\|_{\ell^1_\delta(C^r(X))} \le \eps_\vartheta,
$$
with suitable constants $\eps_\tau$ and $\eps_\vartheta$ when $s=0$.

We call an $L^2_k$ pair $(A,\Phi)$ in the {\em pre-configuration space\/},
$$
\tsC_{W,E} := \sA_E\times L^2_k(X,W^+\otimes E),
$$
a {\em $\PU(2)$ monopole\/} if $\fs(A,\Phi)=0$, where the $\ssG$
equivariant map $\fs:\tsC_{W,E}\to L^2_k(\La^+\otimes\su(E))\oplus
L^2_k(W^-\otimes E)$ is defined by
\begin{equation}
\fs(A,\Phi)
:=
\begin{pmatrix}
F_A^+ - (\id + \tau_0\otimes\id_{\su(E)} + \vectau\cdot\vecfm(A))
(\rho^+)^{-1}(\Phi\otimes\Phi^*)_{00} 
\\
D_A\Phi + \rho(\vartheta_0)\Phi + \vecvartheta\cdot\vecfm(A)\Phi
\end{pmatrix}.
\label{eq:PT}
\end{equation}
We let $M_{W,E} := \fs^{-1}(0)$ be the moduli space of solutions cut out of
the {\em configuration space\/},
$$
\sC_{W,E} := \tsC_{W,E}/\ssG_E,
$$ 
by the section \eqref{eq:PT}, where $u\in\ssG_E$ acts by $u(A,\Phi) :=
(u_*A,u\Phi)$. We note that the perturbations in \eqref{eq:PT} are {\em
zero-order\/}, unlike the first-order perturbations considered
\cite{TelemanMonopole}. They differ from those in
\cite{PTLocal} as they are not `central in $E$',
that is, at some points $x\in X$, they span $\su(E)|_x$.

We let $\sC_{W,E}^{*,0}\subset \sC_{W,E}$ be the subspace of pairs
$[A,\Phi]$ such that $A$ is irreducible and the section $\Phi$ is not
identically zero and set $M_{W,E}^{*,0} = M_{W,E}\cap
\sC_{W,E}^{*,0}$. Note that we have a canonical inclusion
$\sB_E\subset \sC_{W,E}$ given by $[A]\mapsto [A,0]$ and
similarly for the pre-configuration spaces.

The sections $\vectau\cdot\vecfm(A)$ and $\vecvartheta\cdot\vecfm(A)$
vanish at reducible connections $A$ by construction; plainly, the terms in
the $\PU(2)$ monopole equations
\eqref{eq:PT} involving the perturbations $\vectau\cdot\vecfm(A)$ and
$\vecvartheta\cdot\vecfm(A)$ are zero when $\Phi$ is zero. The holonomy
perturbations considered by Donaldson in \cite{DonPoly} are inhomogeneous,
as he uses the perturbations to kill the cokernels of $d_A^+$ directly.
In contrast, the perturbations we consider in \eqref{eq:PT} are homogeneous
and we argue indirectly that the cokernels of the linearization vanish away
from the reducibles and zero-section solutions. 

\section{Uhlenbeck compactness}
\label{sec:UhlenbeckCompactness}
We describe the Uhlenbeck closure of the moduli space of $\PU(2)$ monopoles
and outline the proof of compactness.

\subsection{Statement of main Uhlenbeck compactness results}
\label{subsec:StatementCompactness}
We say that a sequence
of points $[A_\beta,\Phi_\beta]$ in $\sC_{W,E}$ {\em converges\/} to a
point $[A,\Phi,\bx]$ in $\sC_{W,E_{-\ell}}\times\Sym^\ell(X)$, where
$E_{-\ell}$ is a Hermitian two-plane bundle over $X$ such that 
$$
\det(E_{-\ell}) = \det E
\quad\text{and}\quad 
c_2(E_{-\ell}) = c_2(E)-\ell, 
\quad\text{with}\quad 
\ell\in\ZZ_{\ge 0},
$$
if the following hold:
\begin{itemize}
\item There is a sequence of $L^2_{k+1,\loc}$ determinant-one, unitary
bundle isomorphisms $u_\beta:E|_{X\less\bx}\to
E_{-\ell}|_{X\less\bx}$ such that the sequence of monopoles
$u_\beta(A_\beta,\Phi_\beta)$ converges to
$(A,\Phi)$ in $L^2_{k,\loc}$ over $X\less\bx$, and 
\item The sequence of measures  
$|F_{A_\beta}|^2$ converges
in the weak-* topology on measures to $|F_A|^2 +
8\pi^2\sum_{x\in\bx}\delta(x)$.
\end{itemize}

The holonomy-perturbation maps in \eqref{eq:GaugeEquivariantMap}
are continuous with respect to Uhlenbeck
limits, just as are those of \cite{DonPoly}. Suppose $\{A_\beta\}$ is a
sequence in $\sA_E(X)$ which converges to an Uhlenbeck limit $(A,\bx)$ in
$\sA_{E_{-\ell}}(X)\times\Sym^\ell(X)$. The sections
$\vectau\cdot\vecfm(A_\beta)$ and $\vecvartheta\cdot\vecfm(A_\beta)$ then
converge in $L^2_{k+1}(X)$ to a section $\vectau\cdot\vecfm(A,\bx)$ of
$\gl(\La^+)\otimes\so(\su(E_{-\ell}))$ and a section
$\vecvartheta\cdot\vecfm(A,\bx)$ of $\Hom(W^+,W^-)\otimes\fsl(E_{-\ell})$,
respectively.  For each $\ell\ge 0$, the maps of
\eqref{eq:GaugeEquivariantMap} extend continuously to gauge-equivariant
maps
\begin{align}
&\sA_{E_{-\ell}}(X)\times\Sym^\ell(X) 
\to L^2_{k+1}(X,\gl(\La^+)\otimes_\RR\so(\su(E_{-\ell}))), 
\label{eq:GaugeEquivariantExtendedMap}\\
&\sA_{E_{-\ell}}(X)\times\Sym^\ell(X) 
\to L^2_{k+1}(X,\Hom(W^+,W^-)\otimes_\CC\fsl(E_{-\ell})), \notag
\end{align}
given by $(A,\bx)\mapsto \vectau\cdot\vecfm(A,\bx)$ and 
$(A,\bx)\mapsto \vecvartheta\cdot\vecfm(A,\bx)$, 
which are $C^\8$ on each $C^\8$ stratum determined by $\Sym^\ell(X)$.

Our construction of the Uhlenbeck compactification for
$M_{W,E}$ requires us to consider moduli spaces
$$
\bM_{W,E_{-\ell}} \subset \sC_{W,E_{-\ell}}\times\Sym^\ell(X)
$$
of triples $[A,\Phi,\bx]$ given by the zero locus of the $\ssG_{E_{-\ell}}$
equivariant map
$$
\fs:\sC_{W,E_{-\ell}}\times\Sym^\ell(X) \to 
L^2_k(\La^+\otimes\su(E_{-\ell}))\oplus L^2_k(W^-\otimes E_{-\ell})
$$
defined as in \eqref{eq:PT} except using the perturbing sections
$\vectau\cdot\vecfm$ and $\vecvartheta\cdot\vecfm$ in
\eqref{eq:GaugeEquivariantExtendedMap} instead of those in
\eqref{eq:GaugeEquivariantMap}. We call
$\bM_{W,E_{-\ell}}$ a {\em lower-level\/} moduli space if $\ell>0$ and
call $\bM_{W,E_{-0}} = M_{W,E}$ the {\em top\/} or {\em highest level\/}.

In the more familiar case of the unperturbed $\PU(2)$ monopole equations
\eqref{eq:IntroUnpertPT}, the spaces $\bM_{W,E_{-\ell}}$
would simply be products $M_{W,E_{-\ell}}\times \Sym^\ell(X)$.  In general,
though, the spaces $\bM_{W,E_{-\ell}}$ are not products when $\ell > 0$ due
to the slight dependence of the section $\fs(A,\Phi,\bx)$ on the points
$\bx\in\Sym^\ell(X)$ through the perturbations $\vectau\cdot\vecfm$
and $\vecvartheta\cdot\vecfm$. A similar phenomenon is encountered in
\cite[\S 4(iv)--(vi)]{DonPoly} for the case of the extended anti-self-dual
equations, where holonomy perturbations are also employed in order to
achieve transversality.

We define $\barM_{W,E}$ to be the Uhlenbeck closure of $M_{W,E}$ in the
space of ideal $\PU(2)$ monopoles,
$$
IM_{W,E} := \bigcup_{\ell=0}^N \bM_{W,E_{-\ell}} 
\subset 
\bigcup_{\ell=0}^N \left(\sC_{W,E_{-\ell}}\times\Sym^\ell(X)\right)
$$
for any integer $N\ge N_p$, where $N_p$ is a sufficiently large constant to
be specified below.  

\begin{thm}
\label{thm:Compactness}
\cite{FL1}
Let $X$ be a closed, oriented, smooth four-manifold with $C^\8$ Riemannian
metric, \spinc structure $(\rho,W)$ with 
\spinc connection, and a Hermitian two-plane
bundle $E$ with unitary connection on $\det E$.  Then there is a positive
integer $N_p$, depending at most on the curvatures of the fixed connections
on $W$ and $\det E$ together with $c_2(E)$, such that for all $N\ge N_p$
the topological space $\barM_{W,E}$ is second countable, Hausdorff, 
compact, and given by the closure
of $M_{W,E}$ in $\cup_{\ell=0}^{N}\bM_{W,E_{-\ell}}$.
\end{thm}

\begin{rmk}
The existence of an Uhlenbeck compactification for the moduli space of
solutions to the unperturbed $\PU(2)$ monopole equations
\eqref{eq:IntroUnpertPT} was announced by Pidstrigach
\cite{PTCambridge} and an argument was outlined in \cite{PTLocal}.
A similar argument for the equations 
\eqref{eq:IntroUnpertPT} was outlined by Okonek and Teleman in
\cite{OTQuaternion}. Theorem \ref{thm:Compactness} yields the standard
Uhlenbeck compactification for the system \eqref{eq:IntroUnpertPT} and for
the perturbations of \eqref{eq:IntroUnpertPT} described in
\cite{FeehanGenericMetric,TelemanGenericMetric} --- see Remark
\ref{rmk:TransversalityThm}. 
An independent proof of Uhlenbeck compactness for
\eqref{eq:IntroUnpertPT} and certain perturbations of these equations
is given in \cite{TelemanMonopole,TelemanGenericMetric}.
\end{rmk}

\subsection{Outline of the proof of Theorem \ref{thm:Compactness}}
\label{subsec:OutlineCompactness}
The proofs of existence of an Uhlenbeck compactification for the moduli
spaces of solutions to \eqref{eq:IntroUnpertPT} and \eqref{eq:PT}
proceed along similar lines. The common thread is the
use of Bochner-Weitzenb\"ock formulas to give bounds on the `energy' and on the
self-dual component of the curvature 
\begin{align}
\int_X\left(|F_A|^2+|\Phi|^4+|\cov_A\Phi|^2\right)\,dV \le M < \8,
\label{eq:EnergyBound}
\\
\|F_A^+\|_{C^0(X)} \le K,
\label{eq:SelfDualCurvBound}
\end{align}
for any point $[A,\Phi]\in M_{W,E}$, with constant $K$ depending only the
data in the hypotheses of Theorem \ref{thm:Compactness}.
Granted the bounds \eqref{eq:EnergyBound} and \eqref{eq:SelfDualCurvBound},
which are not obvious in the presence of infinite sequences of holonomy
perturbations, the proof of Theorem \ref{thm:Compactness} then follows
the example of \cite[pp. 163--165]{DK}, \cite[Theorem III.3.2]{FrM}
(due to Taubes), and \cite[Proposition 4.4]{TauPath} for the moduli space
of anti-self-dual connections. The scale invariance of the $\PU(2)$
monopole equations is used in
much the same way that the conformal invariance of the anti-self-dual
equation is exploited in \cite{DK}: if $(A,\Phi)$ is a solution to
\eqref{eq:PT} for a Riemannian metric $g$, then $(A,\la\Phi)$ is a solution
for the metric $\la^{-2}g$.

The first step is to prove $C^\8$ regularity of $L^2_1$ solutions to the
$\PU(2)$ monopole equations \eqref{eq:PT} coupled with a Coulomb
gauge-fixing condition with $C^\8$ reference pair when
$(\vectau,\vecvartheta)=(0,0)$ and $C^\8$ regularity of $L^2_2$ solutions
for any $C^\8$ perturbations $(\vectau,\vecvartheta)$.
We next record a generalization of the usual Bochner-Weitzenb\"ock identity
for the Dirac operator on $\Om^0(W^+)$ (see
\cite{KMThom,MorganSWNotes,Witten}) to the Dirac operator on
$\Om^0(W^+\otimes E)$: 
\begin{align}
D_A^*D_A &= \cov_A^*\cov_A + \quarter R
+ \rho^+(F_A^+) + \half\rho^+(F^+(A_{\det E}) + F^+(A_{\det W^+})), 
\label{eq:BWDirac+}\\
D_AD_A^* &= \cov_A^*\cov_A + \quarter R
+ \rho^-(F_A^-) + \half\rho^-(F^+(A_{\det E}) + F^-(A_{\det W^+})), 
\label{eq:BWDirac-}
\end{align}
where $R$ is the scalar curvature of the metric $g$.
As in the case of the $\U(1)$ monopole equations
\cite{KMThom,Witten}, the Bochner-Weitzenb\"ock
identity \eqref{eq:BWDirac+}, integration by parts, the estimates on the
holonomy perturbations given by \eqref{eq:CompactEstPertC0},
and the maximum principle applied to $|\Phi|^2$, yield {\em a priori} estimates
for $\Phi$ and $F_A$ when $(A,\Phi)$ is a solution to 
the $\PU(2)$ monopole equations \eqref{eq:PT}:

\begin{lem}
\label{lem:AprioriEstAPhi}
There is a positive constant $K$, depending only on the data in the
hypotheses of Theorem \ref{thm:Compactness}, such that the following holds.
If $(A,\Phi)$ is an $L^2_1$ solution to \eqref{eq:PT}, then 
\begin{equation}
\|\Phi\|_{L^4(X)}+\|\cov_A\Phi\|_{L^2(X)} \le K
\quad\text{and}\quad
\|F_A\|_{L^2(X)} \le K.
\end{equation}
If $(A,\Phi)$ is a $C^1$ solution to \eqref{eq:PT}, then
\begin{equation}
\|\Phi\|_{C^0(X)} \le K 
\quad\text{and}\quad
\|F_A^+\|_{C^0(X)} \le K.
\end{equation}
\end{lem}

\begin{rmk}
A modification of the vanishing argument in \cite{Witten} shows
that the only solutions to \eqref{eq:PT} on $S^4$, with its standard
metric, are of the form $(A,0)$ with $A$ anti-self-dual.
\end{rmk}

The bounds \eqref{eq:EnergyBound} and \eqref{eq:SelfDualCurvBound} are thus
given by Lemma \ref{lem:AprioriEstAPhi}. 
The energy bound \eqref{eq:EnergyBound} together with the $C^\8$ regularity
of $L^2_1$ solutions to \eqref{eq:PT} and a local Coulomb gauge-fixing
condition are then used to prove removability of point singularities for
$\PU(2)$ monopoles when $(\vectau,\vecvartheta)=0$ following
\cite{UhlRem,UhlChern}. One only needs
removability of singularities in this case as, by construction, the
holonomy perturbations vanish near points of concentrated curvature. The
proof of Theorem \ref{thm:Compactness} is now completed by the usual
measure-theoretic arguments, together with the local Coulomb gauge-fixing
theorem and the patching arguments of \cite{UhlLp}.  One can either use the
Chern-Simons function, as in \cite[\S 4.4]{DK}, to show that the singular
points have integer multiplicities or use the cutting-off procedure of
\cite{TauPath}.

\section{Transversality}
\label{sec:Transversality}
In this section we explain why, for generic perturbation parameters
$(\tau_0,\vectau,\vecvartheta)$, the moduli space $M_{W,E}$ of solutions to
the perturbed $\PU(2)$  monopole equations \eqref{eq:PT} is a smooth
manifold on the complement of the subspaces of zero-section and reducible
solutions.

\subsection{Statement of main transversality results}
\label{subsec:StatementTransversality}
The space $\Sym^\ell(X)$ is smoothly stratified, the strata being
enumerated by partitions of $\ell$. If $\Sigma\subset\Sym^\ell(X)$ is a
smooth stratum, we define
$$
\bM_{W,E_{-\ell}}|_\Sigma
:=\{[A,\Phi,\bx]\in\bM_{W,E_{-\ell}}: \bx\in\Sigma \},
$$
with $\bM_{W,E_{-0}} := M_{W,E}$ when $\ell=0$. 

\begin{thm}
\label{thm:Transversality}
\cite{FL1}
Let $X$ be a closed, oriented, smooth four-manifold with $C^\8$ Riemannian
metric, \spinc structure $(\rho,W)$ with \spinc connection, and a Hermitian
line bundle $\det E$ with unitary connection.  Then there exists a
first-category
subset of the space of $C^\8$ perturbation parameters such that the
following holds: For each $4$-tuple
$(\tau_0,\vartheta_0,\vectau,\vecvartheta)$ in 
the complement of this first-category subset, integer
$\ell\ge 0$, and smooth stratum $\Sigma\subset\Sym^\ell(X)$, the moduli
space $\bM^{*,0}_{W,E_{-\ell}}(\Sigma;\tau_0,\vectau,\vecvartheta)$ is a
smooth manifold of the expected dimension
\begin{align*} 
\dim\bM^{*,0}_{W,E_{-\ell}}|_\Sigma
&= \dim M^{*,0}_{W,E_{-\ell}} + \dim\Sigma 
\\
&= -2p_1(\su(E_{-\ell}))-\threehalf(e(X)+\sigma(X))  + \dim\Sigma\\
&\quad + \half p_1(\su(E_{-\ell}))+\half((c_1(W^+)+c_1(E))^2-\sigma(X))-1,
\end{align*}
where $c_1(E_{-\ell})=c_1(E)$ and $c_2(E_{-\ell})=c_2(E)-\ell$.
\end{thm}

\begin{rmk}
\begin{enumerate}
\item Over the complement of a first-category subset in $\Sigma$, 
comprising the regular values of the projection maps onto the second factor
$\Sigma$, the projection $\bM_{W,E_{-\ell}}^{*,0}|_\Sigma\to\Sigma$ is a
smooth fiber bundle with fibers $\bM_{W,E_{-\ell}}^{*,0}|_{\bx}$, for
$\bx\in\Sigma$.  Indeed, the only tangent vectors in each stratum $\Sigma$
which might not appear in the image of the projection are those arising
from the radial vector on the annuli $4B_j-2\barB_j$.  This observation
shows that the projection from $\bM^{*,0}_{W,E_{-\ell}}|_\Sigma$ to
$\Sigma$ is transverse to certain submanifolds of $\Sigma$ which allows
dimension-counting arguments \cite{FL2}. An approach to dimension counting
in the presence of holonomy perturbations is also discussed by Donaldson in
\cite[pp. 282--287]{DonPoly}.
\item Although not required by Theorem
\ref{thm:Transversality}, a choice of generic 
Riemannian metric on $X$ ensures that the moduli space
$M_E^{\asd,*}$ of irreducible
anti-self-dual connections on $\su(E)$
is smooth and of the expected dimension \cite{DK,FU}, although the points  
of $M_E^{\asd,*}$ need not be regular points of $M^*_{W,E}$ as the
linearization of \eqref{eq:PT} need not be surjective \cite{FL2}.
\item A choice of generic parameter $\tau_0$ ensures that the moduli spaces
$M^{\red,0}_{W,E,L_1}$ of non-zero-section solutions to
\eqref{eq:PT} which are reducible with respect to the splitting
$E=L_1\oplus (\det E)\otimes L_1^*$ are smooth and of the expected
dimension \cite{FL2}, although the points  
of $M^{\red,0}_{W,E,L_1}$ need not be regular points of $M^0_{W,E}$ since the
linearization of \eqref{eq:PT} need not be surjective, as we see in \cite{FL2}.
\end{enumerate}
\end{rmk}

\begin{rmk}
\label{rmk:TransversalityThm}
Different approaches to the question of transversality for the $\PU(2)$
monopole equations \eqref{eq:IntroUnpertPT} with generic perturbation
parameters have been considered by Pidstrigach and Tyurin in \cite{PTLocal}
and by Teleman in \cite{TelemanMonopole}; see \cite{FL1} for further
details.  More recently, a new approach to transversality for
\eqref{eq:IntroUnpertPT} has been discovered independently by the first
author \cite{FeehanGenericMetric} and by Teleman
\cite{TelemanGenericMetric}: the method uses only the perturbations
$(\tau_0,\vartheta_0)$ together with perturbations of the Riemannian metric
on $X$ and compatible Clifford map. Although these perturbations are
simpler than those of Theorem \ref{thm:Transversality}, the proofs of
transversality are lengthy and difficult.
\end{rmk}

While the description of the holonomy perturbations outlined above may
appear fairly complicated at first glance, in practice they do not present
any major difficulties beyond those that would be encountered if simpler
perturbations not involving the bundle $\su(E)$ (such as the Riemannian
metric on $X$ or the connection on $\det W^+$) were sufficient to achieve
transversality \cite{FL2, FL3, FL4}. We note that related transversality
and compactness issues have been recently considered in approaches to
defining Gromov-Witten invariants for general symplectic manifolds
\cite{LiTian,RuanGW,Siebert}.

\subsection{Outline of the proof of Theorem \ref{thm:Transversality}}
\label{subsec:OutlineTransversality}
We shall outline the proof when
$\ell=0$ and indicate the very small modification required when $\ell>0$.
Let $\sP =
\sP_0\oplus\sP_\tau\oplus\sP_\vartheta$ denote our Banach space of $C^r$
perturbation parameters. Define a $\ssG_E$-equivariant map
$$
\fS:\sP\times \tsC_{W,E}
\to L^2_{k-1}(\La^+\otimes\su(E))\oplus L^2_{k-1}(W^-\otimes E)
$$
by setting
\begin{equation}\label{eq:ParamMonopole}
\fS(\tau_0,\vartheta_0,\vectau,\vecvartheta,A,\Phi) 
:= 
\begin{pmatrix}
F_A^+
-(\id+\tau_0\otimes\id_{\su(E)}+
\vectau\cdot\vecfm(A))\rho^{-1}(\Phi\otimes\Phi^*)_{00} \\
D_A\Phi+ \rho(\vartheta_0)\Phi
+\vecvartheta\cdot \vecfm(A)\Phi
\end{pmatrix},
\end{equation}
where $\ssG_E$ acts trivially on the space of perturbations $\sP$, and so
the parametrized moduli space $\fM_{W,E}:=\fS^{-1}(0)/\ssG_E$ is a subset
of $\sP\times\sC_{W,E}$.  We let $\fM_{W,E}^{*,0}:=\fM_{W,E}\cap (\sP\times
\sC_{W,E}^{*,0})$.
The $\ssG_E$-equivariant map $\fS$ defines, in the usual way, a section of
a Banach vector bundle $\ubarfV$ over $\sP\times \sC_{W,E}^{*,0}$,
and a Fredholm section $\fs :=
\fS(\tau_0,\vartheta_0,\vectau,\vecvartheta,\cdot)$ of the Banach vector
bundle $\fV:= 
\ubarfV|_{(\tau_0,\vartheta_0,\vectau,\vecvartheta)}$ over
$\{\tau_0,\vartheta_0,\vectau,\vecvartheta\}\times\sC_{W,E}^{*,0}$ for each
point $(\tau_0,\vartheta_0,\vectau,\vecvartheta)\in\sP$. The principal
goal, then, is to prove

\begin{thm}\label{thm:SmoothParamModuliSpace}
\cite{FL1}
The zero set in $\sP\times \sC_{W,E}^{*,0}$ of the section $\fS$ is regular
and, in particular, the moduli space $\fM^{*,0}_{W,E}$ is a smooth Banach
submanifold of $\sP\times \sC_{W,E}^{*,0}$.
\end{thm}

Accepting Theorem \ref{thm:SmoothParamModuliSpace} for the moment, the
Sard-Smale theorem \cite{Smale} (in the form of Proposition 4.3.11 in
\cite{DK}) implies that the zero sets in $\sC^{*,0}_{W,E}$ of the sections
$\fs_{\tau_0,\vartheta_0,\vectau,\vecvartheta}$ are regular for all $C^r$
perturbations $(\tau_0,\vartheta_0,\vectau,\vecvartheta)$ in the complement
in $\sP$ of a first-category subset and so Theorem
\ref{thm:Transversality} is a standard consequence of Theorem
\ref{thm:SmoothParamModuliSpace} for $C^r$ parameters.
The use of perturbation parameters which are only known to be $C^r$, while
necessary to apply the Sard-Smale theorem, proves inconvenient in
practice. The result can be sharpened, however, so that only $C^\8$
parameters are needed by adapting a similar argument due to Taubes
for the Seiberg-Witten moduli space (see \cite{Salamon}). 
This gives Theorem \ref{thm:Transversality} when
$\ell=0$. The proof of Theorem \ref{thm:Transversality} when $\ell>0$ is
the same as for the case $\ell=0$, except that the $\ssG_E$ equivariant map
$\fS$ is now defined on $\sP\times\tsC_{W,E_{-\ell}}\times\Sigma$, where
$\Sigma$ is a smooth stratum of $\Sym^\ell(X)$.

If $[\tau_0,\vartheta_0,\vectau,\vecvartheta,A,\Phi]$ is a point in
$\fS^{-1}(0)$, elliptic regularity for the equations \eqref{eq:PT} ensures
that the point $[A,\Phi]$ has a $C^r$ representative $(A,\Phi)$.
Similarly, if $(v,\psi)$ is an $L^2_{k-1}$ element of the cokernel of
$D\fS$ at the point $(\tau_0,\vartheta_0,\vectau,\vecvartheta,A,\Phi)$,
elliptic regularity for the Laplacian $D\fS(D\fS)^*$, with $C^{r-1}$
coefficients, implies that $(v,\psi)$ is in $C^{r+1}$. {}From the
definition of the holonomy-perturbation sections $\vectau\cdot\vecfm(A)$
and $\vecvartheta\cdot\vecfm(A)$, the linearization of
\eqref{eq:ParamMonopole}, and some linear algebra we obtain

\begin{prop}
\label{prop:Surjectivity}
\cite{FL1}
Suppose $(\tau_0,\vartheta_0,\vectau,\vecvartheta,A,\Phi)$ is a point in
$\fS^{-1}(0)$ with $A$ irreducible and $\Phi\not\equiv 0$.  If $(v,\psi)$
is in the cokernel of $D\fS$ at
$(\tau_0,\vartheta_0,\vectau,\vecvartheta,A,\Phi)$ then $(v,\psi)\equiv 0$
on each ball $B_j$ whose holonomy sections $\{\fm_{j,l,\alpha}(A)\}$ span
$\su(E)|_{B_j}$.
\end{prop}

In order that the holonomy sections $\{\fm_{j,l,\alpha}(A)\}$ span
$\su(E)|_{B_j}$, we see from the construction of \S
\ref{sec:HolonomyPerturbations} that the connection $A|_{2B_j}$ must be
irreducible. We call $U\subset X$ an {\em admissible open set 
for a connection $A$\/}
if it contains all closed balls $\barB_j$ for which $\beta_j[A]>0$.  (The
supports of all the sections $\fm_{j,l,\alpha}(A)$ are contained in
$\cup_{j=1}^{N_b}\barB_j$ and so any open subset of $X$ containing
$\cup_{j=1}^{N_b}\barB_j$ is admissible.) Thus, in order to bring
Proposition \ref{prop:Surjectivity} to bear on the proof of Theorem
\ref{thm:SmoothParamModuliSpace}, we must know that if $(A,\Phi)$ is an
irreducible, non-zero-section solution to \eqref{eq:PT} on $X$, then the
restriction of $A$ to an admissible open set is irreducible.  This crucial
property is implied by the following unique continuation result for
reducible $\PU(2)$ monopoles:

\begin{thm} 
\label{thm:LocalToGlobalReducible}
\cite{FL1}
If $(A,\Phi)$ is a $C^r$ solution of the perturbed $\PU(2)$  monopole
equations \eqref{eq:PT} 
with $\Phi\not\equiv 0$ over a connected, oriented, smooth four-manifold
$X$ with $C^r$ Riemannian metric and $(A,\Phi)$ is reducible on a non-empty
open subset $U\subset X$
with $\barB_j\subset U$ for all $j$ such that $\beta_j[A]>0$,
then $(A,\Phi)$ is reducible on $X$.
\end{thm}

The corresponding property for anti-self-dual connections is
proved as Lemma 4.3.21 in \cite{DK}.  The proof of Donaldson and Kronheimer
relies on the Agmon-Nirenberg unique continuation theorem for an ordinary
differential inequality on a Hilbert space \cite[Theorem 2]{AN}. We show in
\cite{FL1} that their proof adapts to the case
of the $\PU(2)$ monopole equations
\eqref{eq:IntroUnpertPT} or \eqref{eq:PT}, when the initial open set 
where $(A,\Phi)$ is reducible contains the closed balls $\barB(x_j,R_0)$
supporting holonomy perturbations.
The proof of Theorem \ref{thm:LocalToGlobalReducible} has two main steps.
The first and more difficult one is a local extension result for stabilizers of
pairs which are reducible on a suitable ball:

\begin{prop}
\label{prop:LocalUniqueExtension}
\cite{FL1}
Let $X$ be an oriented, smooth four-manifold with $C^r$ Riemannian metric
and injectivity radius $\varrho=\varrho(x_0)$ at a point $x_0$. Suppose
that $0<r_0<r_1\le\half\varrho$. Let $(A,\Phi)$ be a $C^r$ pair solving the
$\PU(2)$ monopole equations \eqref{eq:PT} on $X$. If $u$ is a $C^{r+1}$
gauge transformation of $E|_{B(x_0,r_0)}$ satisfying $u(A,\Phi)=(A,\Phi)$
on $B(x_0,r_0)$, and if either $B_j\cap B(x_0,r_1)=\emptyset$ or
$B_j\subset B(x_0,r_0)$, for all balls $B_j$ 
for which $\beta_j[A]>0$, then there is an extension of $u$ to a
$C^{r+1}$ gauge transformation $\hat u$ of $E|_{B(x_0,r_1)}$ with $\hat
u(A,\Phi)=(A,\Phi)$ on $B(x_0,r_1)$.
\end{prop}

Proposition \ref{prop:LocalUniqueExtension} is a consequence 
(see \cite{Agmon}) of the
Agmon-Nirenberg unique continuation theorem for an ordinary differential
inequality on a Hilbert space \cite[Theorem 2(ii)]{AN}, together with the
observation that the holonomy perturbations vanish outside the balls $B_j$
with $\beta_j[A]>0$ and on any ball $B_j$ with $\beta_j[A]=0$, so the
$\PU(2)$ monopole equations are local with respect to the pair $(A,\Phi)$
on the open subset of $X$ where the stabilizer is to be extended. The
second and much shorter step is to complete the proof of Theorem
\ref{thm:LocalToGlobalReducible} by showing that the local
stabilizers, produced by repeatedly applying Proposition
\ref{prop:LocalUniqueExtension}, fit together to give a non-trivial
stabilizer and thus 
a reducible pair over the whole manifold $X$.

Granted Proposition~\ref{prop:Surjectivity} and
Theorem~\ref{thm:LocalToGlobalReducible}, the proof of
Theorem~\ref{thm:SmoothParamModuliSpace} is easily concluded.  Let
$(\tau_0,\vartheta_0,\vectau,\vecvartheta,A,\Phi)$ be a $C^r$
representative for a point in $\fM_{W,E}^{*,0}$, so that $A$ is irreducible
and $\Phi\not\equiv 0$, and suppose $(v,\psi)$ is in the cokernel of $D\fS$
at the point $(\tau_0,\vartheta_0,\vectau,\vecvartheta,A,\Phi)$. By
definition of $N_b$ we have $\beta_j[A]>0$ for at least one index
$j\in\{1,\dots,N_b\}$. Since $A$ is irreducible on $X$, Theorem
\ref{thm:LocalToGlobalReducible} implies that $A|_{B(x_{j'},2R_0)}$
must be irreducible for some $j'\in\{1,\dots,N_b\}$ such that
$\beta_{j'}[A]>0$; otherwise, $A|_{B(x_{j'},2R_0)}$ would be reducible
for all $j$ such that $\beta_j[A]>0$ and Theorem
\ref{thm:LocalToGlobalReducible} would imply that $A$ would be
reducible over all of $X$, contradicting our assumption that $A$ is
irreducible. Hence, there is at least one ball $B_{j'}$ whose holonomy
sections $\{\fm_{j',l,\alpha}(A)\}$ span $\su(E)|_{B_{j'}}$.  Now
Proposition \ref{prop:Surjectivity} implies that $(v,\psi)\equiv 0$ on the
set $\barB_J(A)$ of all balls $\barB_j$ for which $\beta_j[A]>0$ and
$A|_{2B_j}$ irreducible. Such elements of the cokernel of $D\fS$ have the
unique continuation property by the Aronszajn-Cordes unique continuation
theorem for second-order elliptic inequalities with a real scalar principal
symbol \cite{Aron}. Indeed, the Laplacian $D\fS(D\fS)^*$ is a differential
operator on $X-\barB_J(A)$ and the coefficients
$\delta\fm_{j,l,\alpha}/\delta A$ containing non-local, integral-operator
terms obtained by differentiating holonomies with respect to $A$ (see
\cite[\S A.2]{FL1} or \cite[\S 8]{TauCasson}) are supported on $B_J(A)$,
precisely where $(v,\psi)\equiv 0$. Thus, $(v,\psi)$ may be viewed as a
$C^{r+1}$ element of the kernel of the second-order, elliptic, differential
operator $D\fS^0(D\fS^0)^*$ with $C^{r-1}$ coefficients over $X$ obtained
by setting all holonomy perturbations and their derivatives equal to zero.
(Without loss of generality, we may scale the Dirac equation in
\eqref{eq:PT} by $1/\sqrt{2}$, so the symbol of $D\fS^0(D\fS^0)^*$ is
$\half$ times the Riemannian metric on $T^*X$.) Hence, $(v,\psi) \equiv 0$
on $X$ and this completes the proof of
Theorem~\ref{thm:SmoothParamModuliSpace}.



\end{document}